\documentclass[11pt]{amsart}

\usepackage{amsthm, amsmath, amssymb,amscd}
\usepackage[all]{xy}
\usepackage{times}
\usepackage{graphicx}
\usepackage{epstopdf}
\usepackage{graphs}

\newtheorem*{theorem}{Theorem}

\title{A counterexample to the maximality of toric varieties}
\author{Valerie Hower}
\begin{document}
\maketitle
\begin{abstract} We present a counterexample to the conjecture of Bihan, Franz, McCrory, and van Hamel concerning the maximality of toric varieties.  There exists a six dimensional projective toric variety $X$ with the sum of the mod $2$ Betti numbers of $X(\mathbb{R})$ strictly less than the sum of the mod $2$ Betti numbers of $X(\mathbb{C})$. 
\end{abstract}
Let $X=X(\mathbb{C})$ be a complex algebraic variety defined by equations with real coefficients.  Complex conjugation $\sigma: X \rightarrow X$ acts on $X$, and the fixed points $X^{\sigma}$ are the real points of $X$, 
denoted $X(\mathbb{R})$.  It is of interest to try to understand the relationship between the topology of $X(\mathbb{C})$ and that of $X(\mathbb{R})$.  If one looks at the mod $2$ Betti numbers, the Smith-Thom inequality states  $$ \sum_i b_i(X(\mathbb{R})) \leq \sum_j b_j(X(\mathbb{C})). $$  We say $X$ is \emph{maximal} or an \emph{M-variety} when equality is achieved.  \\
\indent Toric varieties are defined by equations with integral coefficients.  See Fulton [Ful] for background on toric varieties.  Thus, if $X$ is 
a toric variety, we may consider $X(\mathbb{R})$ the real points of $X$. In [Bih], Bihan et al. define a spectral sequence $G^r_{p,q}$ that converges to the ordinary $\mathbb{Z}_2$ homology of $X(\mathbb{R})$ when $X$ is a projective toric variety.  Here we use different notation and indexing.  We write this spectral sequence  as $\overline{E}^r_{p,q}$, where  $G^r_{p,q}=\overline{E}^r_{p+q,-p}$.
For each $p$ and $q$ we have 
\begin{equation} \overline{E}^1_{p,q} \cong E^2_{p,q} \end{equation} where $E^r_{p,q}$ is the $\mathbb{Z}_2$ Leray-Serre spectral sequence for the moment map of $X(\mathbb{C})$.  
  The spectral sequence $E^r_{p,q}$ converges to the ordinary $\mathbb{Z}_2$ homology of $X(\mathbb{C})$.  
Bihan et al. conjecture that every toric variety is an M-variety.  Their conjecture is based on (1), numerous examples, and the fact that every nonsingular
toric variety is an M-variety.  In this note, we present a counterexample to this conjecture.  The author would like to thank Matthias Franz and Clint McCrory for ongoing discussions and emails.  

\begin{theorem}
There exists a six dimensional projective toric variety which is not an M-variety.
\end{theorem} 
We will define a toric variety from a matroid as in [Gel].
Let $M$ be the rank $3$ matroid $F_7$, the Fano plane.   The affine dependencies of $M$ are depicted below, and the bases of $M$ are the triangles in the diagram. 
  \begin{figure}[htbp]    \begin{center}      \begin{graph}(3.5,3.5)(0,0) 
   \textnode{1}(1.8,3.2){$1$}       
\textnode{7}(1.8,1.07){$7$}
\textnode{2}(1,1.6){$2$}
\textnode{3}(.2,0){$3$}
\textnode{4}(1.8,0){$4$}
\textnode{5}(3.4,0){$5$}
\textnode{6}(2.6,1.6){$6$} 
\edge{1}{2}
\edge{2}{3}
\edge{3}{4}
\edge{4}{5}
\edge{5}{6}
\edge{6}{1}
\edge{1}{7}
\edge{7}{4}
\edge{2}{7}
\edge{7}{5}
\edge{3}{7}
\edge{7}{6} \bow{2}{6}{0.23}\bow{4}{2}{0.27}
\bow{6}{4}{0.27}
   \end{graph}    \end{center}  \end{figure}\\
Consider the projective space $\mathbb{CP}^{27}$ whose coordinates are given by the bases of $M$ $$\{ y_{ijk} \, | \, ijk \in {[7] \choose 3} - \{123,147,156,246, 257,345,367\} \}.$$ We define a linear action of the complex algebraic torus $$T:=(\mathbb{C}^*)^7=\{(t_1,t_2,t_3,t_4,t_5,t_6,t_7)\,| \, t_i\in\mathbb{C}^*\}$$ on  $\mathbb{CP}^{27}$ where the action on the coordinates is 
$$(t_1,t_2,t_3,t_4,t_5,t_6,t_7)\cdot y_{ijk}=t_it_jt_ky_{ijk}.$$  Note that the one dimensional subtorus $\{(t,t,t,t,t,t,t)\,  |\, t\in \mathbb{C}^* \}$ acts trivially.  
By restricting this action to the subtorus $$T^{\prime}:=\{(t_1,t_2,t_3,t_4,t_5,t_6,1)\, |\, t_i\in \mathbb{C}^*\},$$ we obtain an effective action of a six dimensional 
algebraic torus on $\mathbb{CP}^{27}.$  Let $X(\mathbb{C})$ be the closure of the torus orbit $$X(\mathbb{C})=\overline{T^{\prime}\cdot(1,1,1,\ldots ,1)} \subset \mathbb{CP}^{27}.$$  
Then, $X$ is a six dimensional projective toric variety.  The moment polytope for $X$ is the projection of the matroid polytope for $M$
onto the first six coordinates.  
We use $\mathtt{torhom}$ [Fra] to compute the ranks of the entries in $E^2_{p,q}=\overline{E}^1_{p,q}$ for $X $ (2) and the $\mathbb{Z}_2$ Betti numbers for $X(\mathbb{R})$ (3).  
\begin{equation}q
\begin{matrix} & & & & &  &  1\; \\ &&&& &15&0\; \\  &&&&22&0&0 \; \\ &&&6&26&0&0\; \\  &&2&3&9&0&0\;  \\&1&0&1&1&0&0 \; \\ 1 &  0 & 0 & 0 & 0 & 0 & 0 \\  && & p &&& \end{matrix}    \end{equation} 
\begin{equation} \left[\begin{matrix} 1 &1&1&8& 57&15& 1\end{matrix}\right] \end{equation}
As the total grading for the spectral sequence $\overline{E}^r_{p,q}$ is $p$, we see the spectral sequence does not collapse at $\overline{E}^1$.  Moreover, $$\sum_{i=0}^6 b_i(X(\mathbb{R})) = 84.$$
The total grading for the spectral sequence $E^r_{p,q}$ is $p+q$.  Thus, the only possible nonzero higher differential is $d_{4,1}^2:E_{4,1}^2 \rightarrow E_{2,2}^2.$  If $d_{4,1}^2=0$ then the spectral sequence collapses at $E^2$, the $\mathbb{Z}_2$ Betti numbers for $X(\mathbb{C})$ are $$[ 1 \; \; \: 0\; \; \: 1\; \; \: 0\; \; \: 3 \; \; \: 4\; \; \: 15\; \; \: 26\; \; \: 22\; \; \: 0\; \; \:15\; \; \: 0\; \; \: 1]$$
and $$\sum_{i=0}^{12} b_i(X(\mathbb{C})) = 88.$$
If $d_{4,1}^2 \neq 0$ then the $\mathbb{Z}_2$ Betti numbers for $X(\mathbb{C})$ are 
 $$[ 1 \; \; \: 0\; \; \: 1\; \; \: 0\; \; \: 2 \; \; \: 3\; \; \: 15\; \; \: 26\; \; \: 22\; \; \: 0\; \; \:15\; \; \: 0\; \; \: 1]$$

  and $$\sum_{i=0}^{12} b_i(X(\mathbb{C})) = 86.$$
In either case, we obtain $$ \sum_{i=0}^{6} b_i(X(\mathbb{R})) < \sum_{i=0}^{12} b_i(X(\mathbb{C})), $$ and hence $X $ is not an M-variety. 
\\


\begin{thebibliography}{Wag} 
\bibitem[Bih] {Bih}F. Bihan, M. Franz, C. McCrory, J. van Hamel, 
\textit{Is every toric variety an M-variety?}, Manuscripta Math. \textbf{120} (2006), 217-232.
\bibitem[Fra]{Fra} M. Franz, Maple package $\mathtt{torhom},$ version 1.3.0, September 13, 2004, \\ Availible at $\mathtt{http://www-fourier.ujf-grenoble.fr/~franz/maple/torhom.html}$.
\bibitem[Ful] {Ful}W. Fulton,  \textit{Introduction to Toric Varieties}, Princeton University Press, Princeton, NJ, 1993. 
\bibitem[Gel] {Gel}I. M. Gelfand, M. Goresky, R.D. MacPherson, V.V. Segranova, \textit{Combinatorial geometries, convex polyhedra, and Schubert cells}, Adv. Math. \textbf{63} (1988), 301-316. 
\end{thebibliography}
\end{document}